\pdfoutput=1
\RequirePackage{ifpdf}
\ifpdf 
\documentclass[pdftex]{sigma}
\else
\documentclass{sigma}
\fi

\newcommand{\bfa}{\mathbf{a}}
\newcommand{\bfb}{\mathbf{b}}

\newcommand{\bfe}{\mathbf{e}}
\newcommand{\bfg}{\mathbf{g}}
\newcommand{\bfv}{\mathbf{v}}
\newcommand{\bfw}{\mathbf{w}}
\newcommand{\bfx}{\mathbf{x}}
\newcommand{\bfX}{\mathbf{X}}
\newcommand{\cA}{\mathcal{A}}
\newcommand{\cF}{\mathcal{F}}
\newcommand{\cQ}{\mathcal{Q}}
\newcommand{\cT}{\mathcal{T}}
\newcommand{\FF}{\mathbb{F}}

\newcommand{\Ext}{\operatorname{Ext}}
\newcommand{\Gr}{\operatorname{Gr}}
\newcommand{\half}{{\frac{1}{2}}}
\newcommand{\Hom}{\operatorname{Hom}}

\newcommand{\rep}{\operatorname{rep}}
\newcommand{\ZZ}{\mathbb{Z}}

\numberwithin{equation}{section}

\newtheorem{Theorem}{Theorem}
\newtheorem{Corollary}[Theorem]{Corollary}
\newtheorem{Proposition}[Theorem]{Proposition}
 { \theoremstyle{definition}
 }

\begin{document}
\allowdisplaybreaks

\newcommand{\arXivNumber}{1712.08478}

\renewcommand{\thefootnote}{}

\renewcommand{\PaperNumber}{007}

\FirstPageHeading

\ShortArticleName{Some Consequences of Categorification}

\ArticleName{Some Consequences of Categorification\footnote{This paper is a~contribution to the Special Issue on Cluster Algebras. The full collection is available at \href{https://www.emis.de/journals/SIGMA/cluster-algebras.html}{https://www.emis.de/journals/SIGMA/cluster-algebras.html}}}

\Author{Dylan RUPEL~$^\dag$ and Salvatore STELLA~$^\ddag$}

\AuthorNameForHeading{D.~Rupel and S.~Stella}

\Address{$^\dag$~Pasadena Unified School District, Math Academy, Pasadena, CA 91101, USA}
\EmailD{\href{mailto:rupel.dylan@pusd.us}{rupel.dylan@pusd.us}}

\Address{$^\ddag$~Universit\`a degli studi di Roma ``La Sapienza'', Dipartimento di Matematica\\
\hphantom{$^\dag$}~``G. Castelnuovo'', P.le Aldo Moro, 5 -- 00185 Rome, Italy}
\EmailD{\href{mailto:stella@mat.uniroma1.it}{stella@mat.uniroma1.it}}
\URLaddressD{\url{http://www1.mat.uniroma1.it/people/stella/}}

\ArticleDates{Received October 16, 2019, in final form January 21, 2020; Published online January 30, 2020}

\Abstract{Several conjectures on acyclic skew-symmetrizable cluster algebras are proven as direct consequences of their categorification via valued quivers. These include conjectures of Fomin--Zelevinsky, Reading--Speyer, and Reading--Stella related to $\mathbf{d}$-vectors, $\mathbf{g}$-vectors, and $F$-polynomials.}

\Keywords{acyclc cluster algebras; categorification; valued quivers}

\Classification{13F60; 16G20}

\renewcommand{\thefootnote}{\arabic{footnote}}
\setcounter{footnote}{0}

\section{Introduction}

The categorification of skew-symmetric cluster algebras using representations of quivers was initiated by Marsh, Reineke, and Zelevinsky \cite{marsh-reineke-zelevinsky} for Dynkin quivers. With the advent of cluster characters \cite{caldero-chapoton}, the subject has exploded as an industry all its own, leading to the publication of numerous articles including \cite{buan-marsh-reineke-reiten-todorov,caldero-chapoton-schiffler, caldero-keller2,caldero-keller, derksen-weyman-zelevinsky,geiss-leclerc-schroer,palu,plamondon,qin,rupel1,rupel2}, just to name a few.

 The main idea is to understand the combinatorics of cluster mutations in terms of a category of representations of a quiver, and in particular to obtain an interpretation of Laurent expansions of non-initial cluster variables as cluster characters. These are generating functions for certain geometric invariants (e.g., Euler characteristics, point counts over finite fields, Poincar\'e polynomials, etc.) of varieties of subrepresentations in an associated quiver representation.

 Using categorification, many structural conjectures on skew-symmetric cluster algebras, their cluster variables, and the associated combinatorics of mutations have been established. The main goal of this note is to observe that the categorification of acyclic skew-symmetrizable (quantum) cluster algebras \cite{rupel1,rupel2} naturally leads to proofs of many of the same conjectures.

 Alternative approaches to categorification of skew-symmetrizable cluster algebras using quivers with automorphism or group species with potential have been introduced by Demonet \cite{demonet1,demonet2}. Another approach using species with potential has been pursued by Labardini-Fragoso and Zelevinsky~\cite{labardini-zelevinsky}, see Section~15 there for a broader overview of current developments. These formalisms have allowed to prove many conjectures on cluster algebras that we thus omit here even though they should be deducible in our setting.

 Since the appearance of these works, new conjectures on cluster algebras, particularly ones related to denominator vectors, have been made and our goal here is to resolve several of these. In some cases our proofs will follow by direct translation of arguments given in the skew-symmetric case.
 \begin{Theorem} \label{thm:main}
 The following conjectures on cluster algebras are true in the acyclic case: Conjectures~{\rm 4.14(3-4)} from~{\rm \cite{FZ03}}; Conjectures~{\rm 6.11}, {\rm 7.4(1-2)}, {\rm 7.5}, {\rm 7.6}, and~{\rm 7.17} from~{\rm \cite{fomin-zelevinsky4}}; Conjectures~{\rm 3.10} and~{\rm 3.21} from~{\rm \cite{reading-speyer}}; Conjectures~{\rm 2.7}, {\rm 2.8}, and~{\rm 2.9} from~{\rm \cite{reading-stella}}.
 \end{Theorem}

\section{Quantum cluster algebras and quantum cluster characters}\label{sec:qca}

In the interest of brevity, we will assume that the reader is familiar with the basic notions in the theory of cluster algebras and refer to \cite{fomin-zelevinsky4} for background material. Several results are scattered throughout the literature; a compendium focusing on the finite type case can be found in \cite{hohlweg-pilaud-stella}.

Let $B=(b_{ij})$ be a $n\times n$ skew-symmetrizable integer matrix; this means that there exist positive integers $d_i$ such that $d_jb_{ji}=-d_ib_{ij}$ for all $i,j\in[1,n]$. We write $D$ for the diagonal matrix whose~$i$-th diagonal entry is $d_i$ and fix such a matrix once and for all. We assume throughout the paper that $B$ is {\it acyclic}, i.e., up to a simultaneous permutation of rows and columns, $B$ has only non-negative entries above the diagonal. Let $\widetilde{B}$ be the matrix obtained by stacking a $n\times n$ identity matrix below $B$ and write $\cA\big(\bfx,\widetilde{B}\big)$ for the associated cluster algebra; it is a cluster algebra with {\it principal coefficients}.

Since $\widetilde{B}$ has full rank, the algebra $\cA\big(\bfx,\widetilde{B}\big)$ admits a log-canonical Poisson structure compatible with mutations \cite{gekhtman-shapiro-vainshtein} and one may obtain a quantum cluster algebra via a remarkably simple deformation quantization \cite{berenstein-zelevinsky}.

Intuitively, there are only few points to keep in mind. Fix an indetermiate $q$, then the definition of the quantum cluster algebra $\cA_q\big(\bfX,\widetilde B\big)$ can be obtained from that of $\cA\big(\bfx,\widetilde B\big)$ via the following modifications:
 \begin{itemize}\itemsep=0pt
\item The initial cluster $\bfX=(X_i)$ consists of $2n$ variables which $q$-commute, i.e., there exists a~skew-symmetric $2n\times 2n$ matrix $\Lambda=(\lambda_{ij})$ so that $X_iX_j=q^{\lambda_{ij}}X_jX_i$.
\item Let $\cT_\Lambda$ be the quantum torus algebra generated by $X_i^{\pm1}$, $i\in[1,2n]$, over the ring $\ZZ\big[q^{\pm\half}\big]$. It admits an anti-involution (called the {\it bar-involution}) fixing each $X_i$ and interchan\-ging~$q$ with $q^{-1}$. Every non-initial cluster variable should also be ``bar-invariant'' and this uniquely determines the power of $q^\half$ by which to multiply each monomial in the exchange relations. More precisely, using the bar-invariant monomial basis~$X^\bfa$, $\bfa\in\ZZ^{2n}$, of~$\cT_\Lambda$ given by
 \[
 X^\bfa=q^{-\half\sum\limits_{i<j}\lambda_{ij}a_ia_j}X_1^{a_1}\cdots X_{2n}^{a_{2n}},
 \]
we may write the mutated variables as $X'_k=X^{\bfb_+^k-\varepsilon_k}+X^{\bfb_-^k-\varepsilon_k}$, where $\varepsilon_k$ denotes the $k$-th standard basis vector of $\ZZ^{2n}$ and the $k$-th column of $\widetilde{B}$ decomposes as $\bfb^k=\bfb^k_+-\bfb^k_-$ for minimal non-negative vectors $\bfb^k_+,\bfb^k_-\in\ZZ_{\ge0}^{2n}$.
\item Each cluster $\bfX'$ obtainable from $\bfX$ by a sequence of mutations should again generate a~quantum torus, i.e., consist of $q$-commuting variables.
 This forces a compatibility condition between $\widetilde{B}$ and $\Lambda$. That is the matrix $\Lambda$ has to have the following form for some skew-symmetric $n\times n$ matrix $\Lambda_0$:
 \[
 \Lambda = \left[ \begin{matrix}
 \Lambda_0 & -\Lambda_0B-D\\
 -B^T\Lambda_0+D & B^T\Lambda_0B+B^TD
 \end{matrix}
 \right]
 \]
 and thus $\widetilde{B}^T\Lambda=\big[\begin{matrix} D & \boldsymbol{0}\end{matrix}\big]$. An easy calculation shows that this condition naturally reproduces under mutations.
 Note that this compatibility condition is identical to the one required to define a compatible Poisson structure.
 \end{itemize}

 Let $\cF_\Lambda$ denote the skew-field of fractions of $\cT_\Lambda$. The {\it quantum cluster algebra} $\cA_q\big(\bfX,\widetilde{B}\big)$ is the~$\ZZ\big[q^{\pm\half}\big]$-subalgebra of $\cF_\Lambda$ generated by all quantum cluster variables obtainable from the initial cluster $\bfX$ by a sequence of mutations. By work of Berenstein and Zelevinsky \cite{berenstein-zelevinsky}, the famous Laurent phenomenon holds in this adapted setting: $\cA_q\big(\bfX,\widetilde{B}\big)$ is a subalgebra of~$\cT_\Lambda$. Similar to the case of classical cluster algebras, there are formulas to describe these quantum Laurent expansions in terms of the representation theory of valued quivers that we recall below.

To the matrix $B$ is associated a quiver $Q$ with vertices $\{1,\ldots,n\}$ and $\gcd(b_{ij},b_{ji})$ arrows $i\to j$ whenever $b_{ij}<0$; by our assumption this quiver is acyclic. Write $\widetilde{Q}$ for the quiver obtained from $Q$ by attaching to each vertex $i$ an additional vertex $n+i$ via a single arrow $i\to n+i$. Let~$\widetilde{D}=\operatorname{diag}(D,D)$. We will be interested in the representation theories of the valued quivers~$(Q,D)$ and $(\widetilde{Q},\widetilde{D})$.

Let $\FF$ be a finite field and fix an algebraic closure $\overline{\FF}$ of $\FF$. Write $\FF^{\langle d\rangle}$ for the degree $d$ extension of $\FF$ inside $\overline{\FF}$. Note that this provides a canonical identification of $\FF^{\langle d\rangle}$ as a vector space over~$\FF^{\langle g\rangle}$ whenever $g|d$.

A representation $V=(V_i,V_a)$ of the valued quiver $(Q,D)$ consists of an $\FF^{\langle d_i\rangle}$-vector space~$V_i$ for each vertex $i$ and an $\FF^{\langle\gcd(d_i,d_j)\rangle}$-linear map $V_a\colon V_i\to V_j$ for each arrow $a\colon i\to j$. The finite-dimensional representations of $(Q,D)$ form a hereditary abelian category denoted by~$\rep_\FF(Q,D)$. The category $\rep_\FF(Q,D)$ naturally embeds into the category $\rep_\FF\big(\widetilde{Q},\widetilde{D}\big)$ and we will always identify~$\rep_\FF(Q,D)$ as this full subcategory.

Write $\widetilde{\cQ}$ for the Grothendieck group of $\rep_\FF\big(\widetilde{Q},\widetilde{D}\big)$ and let $\cQ\subset\widetilde{\cQ}$ denote the Grothendieck group of $\rep_\FF(Q,D)$. Since $\widetilde{Q}$ is acyclic, there is a natural identification $\widetilde{\cQ}\cong\ZZ^{2n}$ by taking classes $\alpha_i:=[S_i]$ of the vertex-simple representations as a basis. The Euler--Ringel form
 \[
 \langle V,W\rangle:=\dim_\FF\Hom(V,W)-\dim_\FF\Ext(V,W)
 \]
 on pairs of representations $V$ and $W$ only depends on their classes in $\widetilde{\cQ}$.
 More precisely, it may be computed in the basis of vertex-simples by
 \[
 \langle\alpha_i,\alpha_j\rangle
 =
 \begin{cases}
 d_i & \text{if $i=j$},\\
 d_ib_{ij} & \text{if $b_{ij}<0$},\\
 -d_i & \text{if $j=n+i$},\\
 0 & \text{otherwise}.
 \end{cases}
 \]
 For $\bfe\in\cQ$, define ${}^*\bfe=\sum\limits_{i=1}^n\frac{1}{d_i}\langle\alpha_i,\bfe\rangle\alpha_i\in\cQ$ and $\bfe^*=\sum\limits_{j=1}^{2n}\frac{1}{d_j}\langle\bfe,\alpha_j\rangle\alpha_j\in\widetilde{\cQ}$, noting that these do not depend on the choice of symmetrizing matrix $D$.
 Then the {\it quantum cluster character} of a~representation $V\in\rep_\FF(Q,D)$ with $[V]=:\bfv$ is the element of $\cT_\Lambda$ given by
 \[
 X_V := \sum\limits_{\bfe\in\cQ} |\FF|^{-\half\langle\bfe,\bfv-\bfe\rangle} |{\Gr}_\bfe(V) |X^{\widetilde{B}\bfe-{}^*\bfv},
 \]
 where $\Gr_\bfe(V)$ denotes the set of all subrepresentations $E\subset V$ with isomorphism class $[E]=\bfe$. Note that $\widetilde{B}\bfe={}^*\bfe-\bfe^*$ and, for $j\in[1,n]$, we have ${}^*\alpha_j=\alpha_j+\sum\limits_{i=1}^{n}\min(b_{ij},0)\alpha_i$.

 A representation $V\in\rep_\FF(Q,D)$ is {\it rigid} if $\Ext^1(V,V)=0$. If there exists a rigid representation with dimension vector $\bfv\in\cQ$ for the finite field $\FF$, then there exists a rigid representation of~$(Q,D)$ with dimension vector $\bfv$ over any finite field and thus we refer to the dimension vector~$\bfv$ as being {\it rigid} in this case.
 \begin{Theorem}[\cite{rupel1,rupel2}] \label{th:quantum cluster characters} Let $B$ be any acyclic skew-symmetrizable matrix with symmetrizer $D$ and let $(Q,D)$ be the associated valued quiver.
\begin{enumerate}\itemsep=0pt
\item[$(a)$] For each rigid representation $V\in\rep_\FF(Q,D)$, the quantum cluster character $X_V$ is a~quantum cluster monomial of $\cA_{|\FF|}\big(\bfX,\widetilde{B}\big)$. Moreover, every quantum cluster monomial of $\cA_{|\FF|}\big(\bfX,\widetilde{B}\big)$ not involving initial cluster variables arises in this way, the decomposition of~$V$ into indecomposables mirroring the factorization of the quantum cluster monomial into quantum cluster variables.
\item[$(b)$] For all dimension vectors $\bfv,\bfe\in\cQ$ with $\bfv$ rigid, there exists a polynomial $P_{\bfv,\bfe}(q)$ so that for any finite field $\FF$ and rigid representation $V\in\rep_\FF(Q,D)$ of dimension vector $\bfv$, we have~$|{\Gr}_\bfe(V)|=P_{\bfv,\bfe} (|\FF| )$. These polynomials give a ``generic'' quantum cluster character
 \[
 X_\bfv=\sum\limits_{\bfe\in\cQ} q^{-\half\langle\bfe,\bfv-\bfe\rangle}P_{\bfv,\bfe}(q)X^{\widetilde{B}\bfe-{}^*\bfv},
 \]
 which computes the quantum cluster monomials not involving initial quantum cluster va\-riab\-les of $\cA_q\big(\bfX,\widetilde{B}\big)$ with arbitrary parameter~$q$.
 \end{enumerate}
 \end{Theorem}
 An analogous result was obtained for acyclic skew-symmetric quantum cluster algebras by Qin~\cite{qin} using representations of acyclic quivers and with counting polynomials replaced by Poincar\'e polynomials. Setting $q=1$ in the formula for generic quantum cluster characters from Theorem~\ref{th:quantum cluster characters} gives the following corollary from which we will deduce the main results of this note.
 \begin{Corollary} \label{cor:classical cluster characters} All cluster monomials of the cluster algebra $\cA\big(\bfx,\widetilde{B}\big)$ not involving initial cluster variables are computed by the cluster characters
 \[
 x_\bfv = x^{-{}^*\bfv}\sum\limits_{\bfe\in\cQ} P_{\bfv,\bfe}(1)x^{\widetilde{B}\bfe}
 \]
 as $\bfv$ ranges over all rigid dimension vectors in $\cQ$. In particular, the cluster va\-riab\-le $x_\bfv$ has~$\bfg$-vector given by $-{}^*\bfv$ and its $F$-polynomial is $F_\bfv(y)=\sum\limits_{\bfe\in\cQ} P_{\bfv,\bfe}(1)y^\bfe$.
 \end{Corollary}

 The connection between the representation theory of $(Q,D)$ and the cluster algebra is actually much stronger. For a source (resp. sink) vertex $k$ in $Q$, write $\Sigma_k^-\colon \rep_\FF(Q,D)\to\rep_\FF(\mu_kQ,D)$ (resp.\ $\Sigma_k^+\colon \rep_\FF(Q,D)\to\rep_\FF(\mu_kQ,D)$) for the reflection functor as defined in \cite[Section~2]{dlab-ringel}. In what follows, we will drop the adornment and simply write $\Sigma_k$ for both reflection functors, which one to apply should be clear from context. Write $\rep_\FF^{\langle k\rangle}(Q,D)\subset\rep_\FF(Q,D)$ for the full subcategory consisting of representations which contain no summands isomorphic to the simple~$S_k$.
 \begin{Theorem}[\cite{rupel1}] \label{th:reflection functor} Let $k$ be a sink or a source in $Q$ and let $\bfX'$ be the cluster obtained by mutating the initial cluster $\bfX$ in direction~$k$. For any representation $V\in\rep_\FF^{\langle k\rangle}(Q,D)$, we have $X_V=X'_{\Sigma_kV}$, where~$X'_{\Sigma_kV}$ denotes the quantum cluster character of $\Sigma_kV\in\rep_\FF(\mu_kQ,D)$ in the variables $\bfX'$.
 \end{Theorem}

\section{Deducing the conjectures}
In this section, we apply the results of Section~\ref{sec:qca} to deduce the conjectures in Theorem~\ref{thm:main}.

\begin{Proposition}\label{prop:denominators} For any rigid dimension vector $\bfv\in\cQ$, the denominator vector of $x_\bfv$ is $\bfv$.
 \end{Proposition}
\begin{proof} The proof is identical to that of \cite[Section~4, Corollary~2]{caldero-keller} with appropriate modifications in the valued quiver setting, we recall the details here for convenience of the reader.

First note that for any representation $W$ with $[W]=\bfw$, we have $\frac{1}{d_i}\langle W,I_i\rangle=w_i=\frac{1}{d_i}\langle P_i,W\rangle$, where $I_i$ and $P_i$ denote respectively the injective hull and projective cover of the vertex simple $S_i$. Next using the injective resolution of $S_i$ we see that $\langle W,S_i\rangle\le\langle W,I_i\rangle$ while using the projective resolution gives $\langle S_i,W\rangle\le\langle P_i,W\rangle$. Now consider a subrepresentation $E\subset V$ with $[E]=\bfe$, applying the above considerations we see that the $i$-th component of $\bfe^*+{}^*(\bfv-\bfe)$ is bounded by the $i$-th component of $\bfv$:
 \[\frac{1}{d_i}\langle E,S_i\rangle+\frac{1}{d_i}\langle S_i,V/E\rangle\le\frac{1}{d_i}\langle E,I_i\rangle+\frac{1}{d_i}\langle P_i,V/E\rangle=e_i+(v_i-e_i)=v_i.\]

 To finish the proof, for each $i\in[1, n]$ we must exhibit a subrepresentation $E\subset V$ which realizes this bound.
 To construct such a subrepresentation, let $J_i$ be the set of all vertices~$j$ in $Q$ for which there exists a path (possibly trivial) beginning at $i$ and ending at $j$.
 Now set $E_j=V_j$ for $j\in J_i$ and $E_j=0$ for $j\notin J_i$.
 Recall that in the injective coresolution $0\rightarrow S_i\rightarrow I_i\rightarrow I\rightarrow 0$ the injective representation $I$ has nonzero components only at vertices~$j$ which admit a nontrivial path \emph{to} vertex $i$, while in the projective resolution $0\rightarrow P\rightarrow P_i\rightarrow S_i\rightarrow 0$ the projective representation $P$ has nonzero components only at vertices $j$ which admit a nontrivial path \emph{from} vertex $i$.

 Thus $\langle E,S_i\rangle=\langle E,I_i\rangle-\langle E,I\rangle=\langle E,I_i\rangle$ and $\langle S_i,V/E\rangle=\langle P_i,V/E\rangle-\langle P,V/E\rangle=\langle P_i,V/E\rangle$ and therefore
 \[\frac{1}{d_i}\langle E,S_i\rangle+\frac{1}{d_i}\langle S_i,V/E\rangle=\frac{1}{d_i}\langle E,I_i\rangle+\frac{1}{d_i}\langle P_i,V/E\rangle=e_i+(v_i-e_i)=v_i\]
 as desired.
 \end{proof}

 \begin{Corollary}[{\cite[Conjecture 7.17]{fomin-zelevinsky4} for acyclic initial exchange matrices}]
 For a rigid dimension vector $\bfv\in\cQ$, the denominator vector of $x_\bfv$ is the exponent vector of the monomial obtained by tropically evaluating the corresponding $F$-polynomial at $y_1^{-1},\dots,y_n^{-1}$.
 \end{Corollary}
 \begin{proof}
 It is enough to observe that the $F$-polynomial has a unique monomial of maximal degree and this monomial is divisible by all other monomials.
 Indeed, this is the monomial corresponding to the full subrepresentation and hence it has exponent vector $\bfv$ in $F_\bfv$.
 \end{proof}
 This immediately implies \cite[Conjecture 6.11]{fomin-zelevinsky4} using \cite[Proposition 7.16]{fomin-zelevinsky4}.
 From Proposition~\ref{prop:denominators}, we deduce this weakening of \cite[Conjecture 7.4]{fomin-zelevinsky4} in the acyclic case.
 \begin{Corollary}
 \label{cor:sign_coherence}
 Let $\bfv$ be the denominator vector of a non-initial cluster variable $x_\bfv$ in $\cA\big(\bfx,\widetilde{B}\big)$.
 Then:
 \begin{enumerate}\itemsep=0pt
\item[$1)$]all entries in $\bfv$ are non-negative,
\item[$2)$]$v_i=0$ if and only if there is a seed containing both $x_\bfv$ and the initial cluster variable~$x_i$,
\item[$3)$]$v_i$ depends only on $x_\bfv$ and $x_i$ but not on the seed containing~$x_\bfv$.
 \end{enumerate}
 \end{Corollary}
\begin{proof} By Proposition~\ref{prop:denominators}, the denominator vector of any non-initial cluster variable is a dimension vector and therefore it has non-negative entries. Moreover, its coordinates are independent of the seed containing the cluster variable. This proves points~(1) and~(3).

By the categorification construction \cite{rupel2}, seeds in $\cA\big(\bfx,\widetilde{B}\big)$ are in bijection with local tilting representations $T$. The seed corresponding to $T$ contains the $i$-th initial cluster variable precisely when vertex $i$ is not in the support of~$T$. Any rigid indecomposable representation can be completed to a local tilting representation with the same support and this establish part~(2).
\end{proof}

The key difference in the above corollary compared to the original statement of the conjecture is that in point~(3) we do not claim independence from the seed containing~$x_i$. Nonetheless the corollary immediately implies the weaker statements \cite[Conjecture~7.5]{fomin-zelevinsky4} and \cite[Conjecture~2.9]{reading-stella} for acyclic initial exchange matrices.

Let $E=(e_{ij})$ be the \emph{Euler matrix} given by
 \[
 e_{ij} =
 \begin{cases}
 1 & \text{if $i=j$},\\
 \min(b_{ij},0) & \text{if $i\ne j$}.
 \end{cases}
 \]
 Combining Corollary \ref{cor:classical cluster characters} and Proposition \ref{prop:denominators} we immediately get the following.
 \begin{Proposition}[{\cite[Conjecture 3.21]{reading-speyer}}] \label{prop:d to g} Let $\bfv\in\cQ$ be a rigid dimension vector. Then the~$\bfg$-vector of $x_\bfv$ is $-E\bfv$.
 \end{Proposition}
 \begin{proof} Since the $\bfg$-vector of the cluster monomial $x_\bfv$ is given by $-{}^*\bfv$, the result follows immediately from the definitions of $E$ and the operator ${}^*(-)$.
\end{proof}

Replacing $E$ with its piecewise linear counterpart $\nu_c$ defined in \cite{reading-stella-2}, the last result extends to cover all cluster monomials in $\cA\big(\bfx,\widetilde{B}\big)$.

 \begin{Corollary}[{\cite[Conjecture 7.6]{fomin-zelevinsky4}} for acyclic initial exchange matrices] \label{cor:d_vect basis} Different cluster monomials in $\cA\big(\bfx,\widetilde{B}\big)$ have different denominator vectors. Moreover, the denominator vectors of the cluster variables of any given cluster form a $\ZZ$-basis of $\cQ$.
 \end{Corollary}
\begin{proof} Since $B$ is acyclic, by \cite[Remark~7.2 and Proposition~11.6]{demonet2}, different cluster monomials of $\cA\big(\bfx,\widetilde{B}\big)$ have different $\bfg$-vectors.
 Moreover the map $\nu_c$ is invertible over~$\ZZ$~-- in each cluster it is given by a triangular matrix with all the diagonal entries equal to $-1$~-- and the first claim follows from Proposition~\ref{prop:d to g}.

By the same argument, the second claim is a direct consequence of \cite[Proposition~11.5]{demonet2} but can also be deduced as follows. Each seed of $\cA\big(\bfx,\widetilde{B}\big)$ corresponds to a local tilting representation~$T$. By \cite[Lemma~4.3 and Theorem~4.5]{happel-ringel}, the dimension vectors of the summands of~$T$ provide a $\ZZ$-basis for the sublattice on which they are supported. Negative simples corresponding to vertices outside the support complete this basis to a~$\ZZ$-basis of~$\cQ$.
 \end{proof}

\begin{Corollary}[{\cite[Conjecture 2.8]{reading-stella}} for acyclic initial exchange matrices] The mutation of the initial cluster of $\cA\big(\bfx,\widetilde{B}\big)$ at a sink or a source vertex transforms denominator vectors of non-initial cluster variables according to the simple reflection associated to that vertex.
\end{Corollary}
\begin{proof} By Theorem~\ref{th:reflection functor}, the initial cluster mutation at a sink or a source vertex transforms cluster characters according to the associated reflection functor on $\rep_{\FF}(Q,D)$. The result then follows from \cite[Proposition~2.1]{dlab-ringel} and Proposition~\ref{prop:denominators}.
\end{proof}

From this corollary, following \cite[Proposition 2.10]{reading-stella}, we also get \cite[Conjecture 2.7]{reading-stella} for acyclic initial exchange matrices. The following technical result is needed in \cite{rupel-stella-williams}.

 \begin{Proposition} \label{prop:principal F-polynomials} Assume vertex $k$ is a source for $B$. Let $x_\bfv$ be any non-initial cluster variable in $\cA\big(\bfx,\widetilde{B}\big)$ and write $F'_\bfv$ for the tropical evaluation of its $F$-polynomial at $y_j'=y_jy_k^{-b_{kj}-2\delta_{jk}}$. Then~$F'_\bfv=1$ unless $\bfv=[S_k]$ in which case $F'_\bfv=y_k^{-1}$.
 \end{Proposition}
 \begin{proof} Let $V$ be a rigid indecomposable representation with $[V]=\bfv$. The terms of the $F$-polynomial~$F_\bfv(y)$ are labeled by subrepresentations of $[V]$ and the claim is that each of these monomials produces only non-negative exponents when evaluated at the given~$y_j'$.

Since vertex $k$ is a source in $B$, we have $b_{kj}\leq 0$ for $j\in[1,n]$ and therefore only the exponent of~$y_k$ could end up being negative in the tropical evaluation. However, observe that the total exponent of any monomial in $F_\bfv(y')$ is given by applying the simple reflection $s_k$ to the dimension vector of the subrepresentation it corresponds to. In particular, this implies that the only possibility to have a negative exponent is that $S_k$ is a subrepresentation of $V$.

On the other hand, vertex $k$ is a source so that the only rigid indecomposable representation admitting $S_k$ as a subrepresentation is $S_k$ itself and the result follows.
\end{proof}

We turn now our attention towards results dealing with the exchange graph of $\cA\big(\bfx,\widetilde{B}\big)$.

\begin{Proposition}[{\cite[Conjecture 3.10]{reading-speyer}} in the acyclic case] Given a cluster algebra having an acyclic seed, the induced subgraph of its exchange graph consisting of all the seeds containing any fixed collection of cluster variables is connected.
\end{Proposition}
\begin{proof} Since the exchange graph of a cluster algebra with principal coefficients covers the exchange graph of any cluster algebra with the same mutation pattern (cf.~\cite{fomin-zelevinsky4}), it suffices to establish the claim for $\cA\big(\bfx,\widetilde{B}\big)$. The proof given in \cite[Corollary~3]{caldero-keller2} for the weaker statement \cite[Conjecture~4.14(3)]{FZ03} is written in terms of arbitrary hereditary abelian categories. In particular, their proof applies in our situation. More generally, one may replace the rigid indecomposable representation appearing in \cite[Section~5.4]{caldero-keller2} by an arbitrary partial tilting representation (cf.~\cite[Proposition~3]{happel-rickard-schofield}). Then running a similar argument as in the proof of \cite[Theorem~6]{caldero-keller2} proves the general case. These constructions are closely related to the Iyama--Yoshino reduction of triangulated categories (cf.~\cite[Section~4]{iyama-yoshino} or \cite[Section~7.2]{keller}).
 \end{proof}

 \begin{Proposition}[{\cite[Conjecture 4.14(4)]{FZ03}}] The full subgraph of the exchange graph consisting of all seeds having acyclic exchange matrix is connected.
 \end{Proposition}
 \begin{proof} The proof of the analogous statement for skew-symmetric cluster algebras given in \cite[Corollary~4]{caldero-keller2} is also written in terms of arbitrary hereditary abelian categories and thus it applies in our situation as well.
 \end{proof}

\subsection*{Acknowledgements}
 We are grateful to Giovanni Cerulli Irelli for his help with references and to an anonymous referee for pointing out a flaw in the original proof of Corollary~\ref{cor:d_vect basis}. D.R.~was partially supported by an AMS-Simons Travel Grant; S.S.~was supported by ISF grant 1144/16.
D.R.~is grateful to the Department of Physics Mathematics and Astronomy of the California Institute of Technology and to the Department of Mathematics of University of Notre Dame for research support; S.S.~is grateful to the School of Mathematics and Actuarial Science of the University of Leicester, where part of this work was done.

\pdfbookmark[1]{References}{ref}
\LastPageEnding

\end{document}